\documentclass[12pt]{elsarticle}
\usepackage{algorithmic}
\usepackage{amsfonts}
\usepackage{amsmath}
\usepackage{amssymb}
\usepackage{bbm}
\usepackage{color}
\usepackage{eurosym}
\usepackage{graphicx}
\usepackage{latexsym}
\usepackage[T1]{fontenc}
\usepackage{upgreek}
\usepackage{mathtools}

\usepackage{matlab-prettifier}

\newtheorem{theorem}{Theorem}

\newtheorem{problem}{Problem}
\newtheorem{example}{Example}

\newtheorem{remark}{Remark}

\newtheorem{algorithm}{Algorithm}

\newcommand{\bfx}{\mbox{$\mbox{\boldmath $x$}$}} 
\newcommand{\bfy}{\mbox{$\mbox{\boldmath $y$}$}} 

\newcommand{\bfzero}{\mbox{$\mbox{\boldmath $0$}$}} 

\newcommand{\bigzero}{\scalebox{0.9}{\mbox{\large $0$}}}

\newcommand{\rank}{\mathrm{rank\;}}

\newcommand{\ct}{\mbox{$\mathbb C$}}
\newcommand{\oo}{\mbox{$\mathcal O$}}

\def\diag{\mathop{{\rm diag}}\nolimits}

 \def\diag{\mathop{{\rm diag}}\nolimits}

\newcommand*\xbar[1]{%
  \hbox{%
    \vbox{%
      \hrule height 0.5pt 
      \kern0.5ex
      \hbox{%
        \kern-0.1em
        \ensuremath{#1}%
        \kern-0.1em
      }%
    }%
  }%
}

\date{}

\journal{ArXiv: published in Electronic Journal of Linear Algebra}

\begin{document}

\begin{frontmatter}

\title{Optimal approximation of a large matrix by a sum of projected linear mappings on prescribed subspaces}

\author[scg,ph]{Phil Howlett}
\author[scg,mk]{Anatoli Torokhti\corref{cor1}}

\address[scg]{Scheduling and Control Group, Industrial AI Research Centre, School of Mathematics and Information Technology, UniSA STEM,\\ University of South Australia, Mawson Lakes 5095, Australia}

\address[stem]{Retired from School of Mathematics and Information Technology, UniSA STEM,\\ University of South Australia, Mawson Lakes 5095, Australia}

\address[ph]{email:~phil.howlett@unisa.edu.au. url:~http://orcid.org/0000-0003-2382-8137}
\address[mk]{email:~anatoli.torokhti@gmail.com. url:~https://orcid.org/0000-0002-8337-4678}

\cortext[cor1] {Corresponding author}

\begin{abstract}
\label{abs}\noindent { We propose and justify a matrix reduction method} for calculating the optimal approximation of an observed matrix $A \in \ct^{m \times n}$ by a sum $\sum_{i=1}^p \sum_{j=1}^q B_iX_{ij}C_j$ of matrix products where each $B_i \in \ct^{m \times g_i}$ and $C_j \in \ct^{h_j \times n}$ is known and where the unknown matrix kernels $X_{ij}$ are determined by minimizing the Frobenius norm of the error.  The sum can be represented as a bounded linear mapping $BXC$ with unknown kernel $X$ from a prescribed subspace ${\mathcal T} \subseteq \ct^n$ onto a prescribed subspace ${\mathcal S} \subseteq \ct^m$ defined respectively by the collective domains and ranges of the given matrices $C_1,\ldots,C_q$ and $B_1,\ldots,B_p$.  We show that the optimal kernel is $X = B^{\dag}AC^{\dag}$ and that the optimal approximation $BB^{\dag}AC^{\dag}C$ is the projection of the observed mapping $A$ onto a mapping from ${\mathcal T}$ to ${\mathcal S}$.  If $A$ is large $B$ and $C$ may also be large and direct calculation of $B^{\dag}$ and $C^{\dag}$ becomes unwieldy and inefficient. { The proposed method avoids} this difficulty by reducing the solution process to finding the pseudo-inverses { of } a collection of much smaller matrices. This significantly reduces the computational burden.
\end{abstract}

\end{frontmatter}



\section{Introduction}
\label{intro}

Let $A \in \ct^{m \times n}$ be an observed matrix and let
$$
{\mathcal S} = B_1({\mathbb C}^{g_1}) \cup \cdots \cup B_p({\mathbb C}^{g_p}) \subseteq \ct^m \quad \mbox{and} \quad {\mathcal T} = C_1^*({\mathbb C}^{h_1}) \cup \cdots \cup C_q^*({\mathbb C}^{h_q}) \subseteq \ct^n
$$
be prescribed subspaces defined by known matrices $B_i \in \ct^{m \times g_i}$ with $g_i \in {\mathbb N}$ for each $i=1,\ldots,p$ and $C_j \in {\mathbb C}^{h_j \times n}$ with $h_j \in {\mathbb N}$ for each $j=1,\ldots,q$.   In general we assume that $m, n \in {\mathbb N}$ are large and that $g_1,\ldots,g_p \in {\mathbb N}$ and $h_1,\ldots,h_q \in {\mathbb N}$ are small.  Consider the following problem.

\begin{problem}
\label{p1}
Find $X_{ij} \in \ct^{g_i \times h_j}$ for each $i=1,\ldots,p$ and $j=1,\ldots,q$ to solve
\begin{equation}
\label{eq:1.1}
\min_{\substack{X_{ij}\\i=1,\ldots,p, j=1,\ldots,q}} \left \| A - \mbox{$ \sum_{i=1}^p \sum_{j=1}^q$} B_iX_{ij}C_j \right \|_F
\end{equation}
where $\| \cdot \|_F$ denotes the Frobenius norm. $\hfill \Box$
\end{problem}

{  Let $g=g_1+\ldots +g_p$ and $h+h_1 + \ldots + h_q$. }  If we define $B \in {\mathbb C}^{m \times g}$, $C \in {\mathbb C}^{h \times n}$ and $X \in {\mathbb C}^{g \times h}$ by setting
$$
B = \left[ \begin{array}{ccc}
B_1 & \cdots & B_p \end{array} \right], \hspace{0.25cm}C = \left[ \begin{array}{c}
C_1 \\
\vdots \\
C_q \end{array} \right], \hspace{0.25cm} X = \left[ \begin{array}{ccc}
X_{11} & \cdots & X_{1q} \\
\vdots & \ddots & \vdots \\
X_{p1} & \cdots & X_{pq} \end{array} \right]
$$
then Problem \ref{p1} can be stated more succinctly as follows.

\begin{problem}
\label{p2}
Find $X \in \ct^{g \times h}$ to solve
\begin{equation}
\label{eq:1.2}
{\min_{X} \left \| A - BXC \right \|_F}
\end{equation}
where $\| \cdot \|_F$ denotes the Frobenius norm. $\hfill \Box$
\end{problem}

\subsection{Preliminaries}\label{ss:fpd}

It is well known \cite{fri1} that $X \in \ct^{g \times h}$ is a solution to (\ref{eq:1.2}) if and only if it is also a solution to
\begin{equation}\label{eq:1.3}
B^*BXCC^* = B^*AC^*.
\end{equation}
Equation (\ref{eq:1.3}) always has a solution but in general the solution is not unique.  If we write $B^{\dag}$ and $C^{\dag}$ for the respective Moore\textendash Penrose inverse matrices then the general solution is given by
\begin{equation}\label{eq:1.4}
X = B^{\dag}AC^{\dag} + (I - B^{\dag}B)M(I - CC^{\dag})
\end{equation}
where $M \in \ct^{g \times h}$ is arbitrary.  For every solution $X$ we have
\begin{equation}\label{eq:1.5}
BXC = BB^{\dag}AC^{\dag}C.
\end{equation}
Therefore the best approximation $BXC$ is uniquely determined.  Since $BB^{\dag}$ is the projection of the space $\ct^m$ onto the subspace ${\mathcal S} \subseteq \ct^m$ and $C^{\dag}C$ is the projection of the space $\ct^n$ onto the subspace ${\mathcal T} \subseteq \ct^n$ we can see that the optimal approximation $BXC \in {\mathcal B}({\mathcal T},{\mathcal S})$ is the projection of the bounded linear mapping $A \in {\mathcal B}(\ct^n,\ct^m)$ onto a bounded linear mapping from ${\mathcal T}$ to ${\mathcal S}$.   Note that if $\bfx \in {\mathcal T}$ and $\bfy = A\bfx \in {\mathcal S}$ then
$$
BXC \bfx = B(B^{\dag}AC^{\dag})C \bfx = BB^{\dag}A(C^{\dag}C \bfx) = BB^{\dag}A \bfx = BB^{\dag} \bfy = \bfy.
$$
In this case the approximate mapping is exact.

 \subsection{ Motivation}\label{s:m}

Our matrix approximation problem is motivated by tasks arising in data processing.  Suppose it is necessary to repeatedly transmit large data matrices from a remote location to a central storage facility. Although there may be numerous portable devices in the field that can transmit data it may also be the case that no single device has the capacity to transmit the complete data set in the time allowed. { This scenario is considered, for example, in \cite{tor22}, \cite{tor23} and the references therein. } More specifically suppose we wish to transmit a large data matrix representing a digital message or an image, from one location to another using an array of transmitters. We assume that each individual transmitter can be used to transmit a small matrix $X_{ij} \in \ct^{h_i \times g_j}$ for each $i=1,\ldots,p$ and $j=1,\ldots,q$.  If $B_i \in \ct^{m \times g_i}$ and $C_j \in \ct^{h_j \times n}$ are known matrices to both sender and receiver, the sender can find each $X_{ij}$ by solving Problem~\ref{p1} and the receiver can reconstruct $A$ approximately from the data array $X = [X_{ij}]$ according to the formula $A \approx \sum_{i=1}^p\sum_{j=1}^q B_iX_{ij}C_j$.  We could regard $X = [X_{ij}]$ as a codified form of the intended message $A$ in which case the matrices $B = [B_i]$ and $C = [C_j]$ would represent decoders.

Alternatively we could imagine that $A$ represents the kernel of a large linear system and that we wish to restrict our analysis to certain key input\textendash output relationships. { See, e.g.,  \cite{tor24} and the references therein.}  Suppose we are interested in those inputs that lie in the space ${\mathcal T}$ defined by the collective domains of the matrices $C_j$ for $j=1,\ldots,q$ and those outputs that lie in the space ${\mathcal S}$ defined by the collective ranges of the matrices $B_i$ for $i=1,\ldots,p$.  The matrix $X$ represents the kernel of an economized linear system $BXC$ which is restricted to the designated input\textendash output relationships.

\subsection{Aim}\label{ncoa}

Our aim is to replace the direct calculation described in (\ref{eq:1.4}) by a more efficient solution procedure and hence reduce the computational burden associated with the direct calculation of generalized inverses of large matrices.

The minimum norm solution to both (\ref{eq:1.2}) and (\ref{eq:1.3}) is uniquely defined by the formula $X_0 = B^{\dag}AC^{\dag}$ where $B \in \ct^{m \times g}$ and $C \in \ct^{h \times n}$.   We have $B^{\dag} = \{B^*B\}^{\dag}B^*$ where $B^*B \in \ct^{g \times g}$ and $C^{\dag} = C^* \{CC^*\}^{\dag}$ where $CC^* \in \ct^{h \times h}$.  Direct calculation of $B^{\dag}$ and $C^{\dag}$ using a standard singular value decomposition algorithm will take respectively $\oo(g^3)$ and $\oo(h^3)$ floating point operations (flops).  When $g, h$ are large this direct calculation becomes unwieldy and inefficient.   We propose { a more efficient alternative} to this direct calculation.  We are primarily interested in the reconstruction of $A$.  In this regard it is important to remember that $BXC = BX_0C$ for all $X$ given by (\ref{eq:1.4}).  Thus it is not necessary to compute the minimum norm solution.

If the subspaces ${\mathcal S}_1 = B_1(\ct^m),\ldots,{\mathcal S}_p = B_p(\ct^m)$ are mutually orthogonal and the subspaces ${\mathcal T}_1 = C_1^*(\ct^n),\ldots,{\mathcal T}_q = C_q^*(\ct^n)$ are also mutually orthogonal then the optimal matrices $X_{ij}$ are determined independently with $X_{ij} = B_i^{\dag}AC_j^{\dag}$ for each $i=1,\ldots,p$ and $j = 1,\ldots,q$.  This is not true in general.

\section{\hspace*{-1mm}Contribution}\label{s:cont}

\hspace*{-1mm}We propose  a new method for numerical solution of the matrix equation $B^*BXCC^* = B^*AC^*$  called the elementary block operations scheme (EBOS).  The EBOS algorithm is described in Section \ref{ss:ebosa} and the corresponding M{\sc atlab} code in presented in \ref{s:mc}.  To begin we write the equation as two separate systems.   We define $Y$ as the solution to the equation $YCC^* = AC^*$ and define $X$ as the solution to the equation $B^*BX = B^*Y$.  The scheme is designed to greatly reduce the computational burden.

$\bullet$ \textbf{Elementary block operations scheme (EBOS).}  We use a result established by Baksalary and Baksalary \cite[Theorem 3, pp 21\textendash 22]{bak1} to justify nonsingular transformations $E = E_1 \cdots E_{p-1}$ and $F = F_{q-1} \cdots F_1$ defined respectively by a product of nonsingular elementary block upper and lower triangular matrices so that the transformed matrices
$$
B^{(p-1)} = BE = \left[ \begin{array}{cccc}
B_1 & B_2^{(1)} & \cdots & B_p^{(p-1)} \end{array} \right] \quad \mbox{and} \quad C^{(q-1)} = FC = \left[ \begin{array}{c}
C_1 \\
C_2^{(1)} \\
\vdots \\
C_q^{(q-1)} \end{array} \right]
$$
have the special property that the Moore\textendash Penrose inverses are given by
$$
[B^{(p-1)}]^{\dag} = \left[ \begin{array}{c}
\rule{0cm}{0.5cm} B_1^{\dag} \\
\rule{0cm}{0.5cm} [B_2^{(1)}]^{\dag} \\
\vdots \\
\rule{0cm}{0.5cm} [B_p^{(p-1)}]^{\dag} \end{array} \right] \quad \mbox{and} \quad [C^{(q-1)}]^{\dag} = \left[ \begin{array}{cccc}
C_1^{\dag} & [C_2^{(1)}]^{\dag} & \cdots & [C_q^{(q-1)}]^{\dag} \end{array} \right].
$$
The reduction also ensures that the matrices $D_B = [B^{(p-1)}]^*B^{(p-1)}$ and $D_C = C^{(q-1)}[C^{(q-1)}]^*$ are block diagonal and that the systems $B^*BX = B^*Y \iff D_B E^{-1}X = [B^{(p-1)}]^*Y$ and $Y^*CC^* = AC^* \iff Y^*F^{-1} D_C = A[C^{(q-1)}]^*$ can be solved by calculating the Moore\textendash Penrose inverses of the block diagonal elements of the transformed coefficients.

\subsection{Benefit}
We will show that {  EBOS} can be easily implemented using the standard general purpose M{\sc atlab} algorithms.  In practice we would expect these general purpose algorithms to be replaced in { the proposed} scheme by more efficient special numerical techniques considered in a number of works.  {We cite \cite{gil1} as a particular instance and refer to \cite{gol4} for a general discussion of the relevant numerical methods.}  The problem of calculating the Moore\textendash Penrose inverses of a $g \times g$ matrix and an $h \times h$ matrix is reduced to a much smaller problem of calculating the Moore\textendash Penrose inverses of a $g_i \times g_i$ matrix for each $i=1,\ldots,p$ and an $h_j \times h_j$ matrix for each $j=1,\ldots,q$.  Thus the number of flops needed for the numerical calculation of the respective Moore\textendash Penrose inverses has been reduced from $\oo(g^3) + \oo(h^3)$ to $\oo(\sum_{i=1}^p g_i^3) + \oo(\sum_{j=1}^q h_j^3)$.
 We make a general comparison of our method with the direct Matlab method by counting the number of flops in a full range of different circumstances.

In problems such as those in biology, ecology, finance, sociology and medicine (see e.g.,  \cite{gene12,lsa10,pca176}), $m,n=O(10^4)$ and greater. For example, measurements of  gene expression can contain tens of hundreds of samples. Each sample can contain tens of thousands of genes. As a result, in such cases the associated  matrices  are very large. In particular, the  well-known phenomenon of the ``curse of dimensionality'' \cite{Bell2003} states that many problems become exponentially difficult in high dimensions.

\section{Literature review}
\label{s:lr}

The topic of matrix approximation has a long history with contributions from many authors.  The motivations for individual works are wide-ranging and the inter-connections are complex.  The early work in the period 1905\textendash 1912 on integral equations \cite{schm, wey1}  was motivated by the development of a systematic theory for solution of integral equations.  The work in  \cite{schm, wey1} essentially established a singular value decomposition (SVD) for bounded integral operators on the infinite-dimensional space of square integrable functions and proved that the best $k$-dimensional approximation to the integral operator was the projected integral operator defined by the characteristic vectors associated with the $k$ largest singular values.  The next contributions in the period 1933\textendash 1941 were motivated by applications to psychometrics \cite{eck1, eck2, hot1, hot2, hou1, scho, wey1, you1, you2} where an observed test score matrix was approximated by a lower rank matrix.  The underlying idea was that linear dependencies in the true score were obscured in the observed score by the addition of random errors.  The most notable results were the Eckart\textendash Young theorem \cite{eck1} which was essentially the Schmidt approximation theorem \cite{schm} in matrix notation and the principal component analysis of Hotelling \cite{hot1, hot2}.  See also later work by Green \cite{gre1} and Sch\"{o}nemann  \cite{scho} which was also motivated by Psychometric considerations.  The SVD was the mathematical foundation for all of this work.  The next results were motivated by fundamental mathematical considerations \cite{mir1} and by development of effective algorithms for numerical calculation of the SVD and the closely related QR decomposition \cite{fra1, fra2, gol1, gol2, kub1}.  A detailed description of the above work can be found in the 1993 survey paper by G W Stewart \cite{ste1}.

Development of a comprehensive theory of generalized inverse matrices was closely related to the conceptual development of the SVD.  The {\em so-called} Moore\textendash  Penrose matrix inverse was first discovered and justified by E.H. Moore \cite{ben1, moo1} and later rediscovered and elaborated by R. Penrose \cite{ben2, pen1}.  See \cite{gol4, lue1} for further discussion.

Friedland and Torokhti \cite{fri1} solved the generalized Eckart\textendash Young problem
\begin{equation}
\label{eq:3.1}
\min_{X \in {\mathcal R}(m,n,k)} \|A -  BXC \|
\end{equation}
where $A \in {\mathbb C}^{m \times n}$, $B \in {\mathbb C}^{m \times g}$ and $C \in {\mathbb C}^{n \times h}$ are known and ${\mathcal R}(m,n,k)$ is the set of all $X \in \ct^{m \times n}$ with $\rank(X) \leq k$.  Liu et al.~\cite{liu1} solved a similar problem for Hermitian or skew-Hermitian matrices.  { See also \cite{wan1,pablo2023}   for  similar forms of the rank constrained matrix approximation.}  Other forms of constrained matrix approximation include the so-called skeleton approximations \cite{bou1, chi1, dri1, gor1}.

The work on matrix approximation has also found applications in signal processing where Hua and Liu \cite{hua1} applied the Eckart\textendash Young theorem to the dimensionality reduction of transmitted random signals.  A similar approximation was established for estimation of random signals in infinite dimensional Hilbert space in \cite{fom1, how1} provided the generalized inverse auto-covariance operator is a bounded operator.  This result was recently further extended by Howlett and Torokhti \cite{how2} to the case where the generalized inverse auto-covariance operator is closed but unbounded.

\section{The elementary block operations scheme (EBOS)}
\label{s:ebos}

EBOS uses a sequence of elementary block row operations to reduce each of the separated systems $B^*BX = B^*Y$ and $YCC^* = AC^*$ to block diagonal form.

\subsection{A theoretical basis for EBOS}
\label{ss:tbe}

We use a key result proved by Baksalary and Baksalary \cite{bak1} in 2007.  See also \cite{bak2} for further discussion.  It is convenient to state the key result in the following form.

\begin{theorem}\label{t:4.1}
Let
$$
R = \left[ \begin{array}{c}
R_1 \\ \hline
R_2 \end{array} \right] \in \ct^{m \times n}
$$
where $R_1 \in \ct^{m_1 \times n}$, $R_2 \in \ct^{m_2 \times n}$ and $m = m_1 + m_2$.  Suppose
$$
S = \left[ \begin{array}{c|c}
R_1^{\dag} & \rule{0cm}{0.5cm} R_2^{\dag} \end{array} \right]
$$
where $R_j^{\dag}$ denotes the Moore\textendash Penrose inverse of $R_j$ for each $j=1,2$.  The following statements are equivalent
\begin{enumerate}
\item $S = R^{\dag}$,
\item $RSR = R$,
\item $R_2^{-1}( \{ \bfzero\})^{\perp} \subseteq R_1^{-1}(\{ \bfzero\}) \iff R_2^*({\mathbb C}^n) \subseteq R_1^{-1}(\{ \bfzero\})$,
\item $R_1^{-1}( \{ \bfzero\})^{\perp} \subseteq R_2^{-1}(\{ \bfzero\}) \iff R_1^*({\mathbb C}^n) \subseteq R_2^{-1}(\{ \bfzero\})$,
\end{enumerate}
where $R_j^{-1}(\{ \bfzero\}) \subseteq \ct^n$ denotes the null space of $R_j$, $R_j^{-1}(\{ \bfzero\})^{\perp} \subseteq \ct^n$ denotes the orthogonal complement of the null space of $R_j$ and $R_j({\mathbb C}^n) \subseteq {\mathbb C}^m$ denotes the range space of $R_j$ in ${\mathbb C}^{m_j}$ considered as a subspace of ${\mathbb C}^m$ for each $j=1,2$. $\hfill \Box$
\end{theorem}

Suppose that $R$ does not satisfy the conditions of Theorem~\ref{t:4.1}.  We can use elementary row operations to transform the matrix $R$ into a matrix $R^{(1)}$ that does satisfy the conditions of Theorem~\ref{t:4.1}.  Define
$$
R^{(1)} = F_1 R = \left[ \begin{array}{c|c}
I & \bigzero \\ \hline
\rule{0cm}{0.5cm} -R_2R_1^{\dag} & I \end{array} \right] \left[ \begin{array}{c}
R_1 \\ \hline
R_2 \end{array} \right] = \left[ \begin{array}{c}
R_1 \\ \hline
\rule{0cm}{0.5cm} R_2(I - R_1^{\dag}R_1)\end{array} \right] = \left[ \begin{array}{c}
R_1^{(1)} \\ \hline
\rule{0cm}{0.5cm} R_2^{(1)} \end{array} \right].
$$
Now $R_2^{(1)} [R_1^{(1)}]^{\dag} = R_2 (I - R_1^{\dag}R_1)R_1^{\dag} = \bigzero$ implies $[R_1^{(1)}]^{\dag}({\mathbb C}^m) \subseteq [R_2^{(1)}]^{-1}(\{ \bfzero\}) \iff R_1^{-1}(\{ \bfzero\})^{\perp} \subseteq R_2^{-1}(\{ \bfzero\})$.  Therefore $
[R^{(1)}]^{\dag} = \left[ \begin{array}{c|c}[R_1^{(1)}]^{\dag} & [R_2^{(1)}]^{\dag} \end{array} \right]$.  For more information see Luenberger \cite[p 164, Figure 6.6]{lue1}.

\begin{example}
\label{ex:1}
Let
$$
R = \left[ \begin{array}{c}
R_1 \\ \hline
R_2 \end{array} \right] = \left[ \begin{array}{rrrrrr}
     1  &    0  &    1   &   0  &    1  &    0\\
    -1  &    0  &    1  &   -1  &    0  &    1\\
     1  &    0  &   -1   &   0  &    1   &   0\\
 \hline
     1  &    0  &    1 &    -1 &     2  &    1\\
     1  &    1  &    1  &    0  &    0  &    0\\
     1  &    2  &   -2  &   -2  &    0 &     1\\
\end{array} \right].
$$
Direct M{\sc atlab} calculations show that $R^{\dag} \neq \left[R_1^{\dag} \mid R_2^{\dag}\right]$.  If we define $R_2^{(1)}= R_2(I - R_1^{\dag}R_1)$ then  M{\sc atlab} returns
$$
R_2^{(1)}  = \left[ \begin{array}{rrrrrr}
    0.0000     &      0   &   0.0000    &       0   &   0.0000    &       0\\
    0.4000  &    1.0000  &    0.0000  &   -0.2000  &   -0.4000  &    0.2000\\
    1.0000  &    2.0000  &   -0.0000  &   -1.0000  &   -1.0000     &      0\\
 \end{array} \right].
$$
Therefore,
$$
R^{(1)}=\left[ \begin{array}{c}
R_1^{(1)} \\ \hline
\rule{0cm}{0.5cm} R_2^{(1)} \end{array} \right]=\left[ \begin{array}{rrrrrr}
    1.0000    &      0   &  1.0000    &      0  &   1.0000     &     0\\
   -1.0000     &     0   &  1.0000  &  -1.0000    &      0  &   1.0000\\
    1.0000    &      0  &  -1.0000   &       0  &   1.0000    &      0\\
    \hline
    0.0000     &     0   &  0.0000    &      0   &  0.0000    &      0\\
    0.4000   &  1.0000   &  0.0000  &  -0.2000  &  -0.4000  &   0.2000\\
    1.0000  &   2.0000  &  -0.0000  &  -1.0000  &  -1.0000    &      0\\
 \end{array} \right].
$$
Now by M{\sc atlab},
$$
 [R_1^{(1)}]^{\dag} = \left[ \begin{array}{rrr}
    0.3000  &   -0.2000   &   0.1000\\
    0.0000  &   -0.0000  &   -0.0000\\
    0.5000  &   -0.0000  &   -0.5000\\
    0.1000  &   -0.4000  &   -0.3000\\
    0.2000  &    0.2000   &   0.4000\\
   -0.1000   &   0.4000   &   0.3000\\
\end{array} \right], \quad
[R_2^{(1)}]^{\dag} = \left[ \begin{array}{rrr}
    0.0000 &   -0.2500 &    0.2500\\
   -0.0000 &    1.2500 &   -0.2500\\
   -0.0000  &   0.0000 &   -0.0000\\
  0.0000 &    2.0000  &  -1.0000\\
   0.0000  &   0.2500  &  -0.2500\\
  0.0000  &   1.7500 &   -0.7500\\
\end{array} \right],
$$
 and
$$
[R^{(1)}]^{\dag} = \left[ \begin{array}{rrr|rrr}
    0.3000  &  -0.2000  &   0.1000 &    0.0000 &   -0.2500  &   0.2500\\
   -0.0000 &   -0.0000 &   -0.0000  &   0.0000 &    1.2500  &  -0.2500\\
    0.5000  &  -0.0000  &  -0.5000  &   0.0000  &   0.0000  &   0.0000\\
    0.1000 &   -0.4000  &  -0.3000 &    0.0000 &    2.0000  &  -1.0000\\
    0.2000 &    0.2000  &   0.4000  &   0.0000  &   0.2500  &  -0.2500\\
   -0.1000 &    0.4000  &   0.3000  &   0.0000  &   1.7500 &   -0.7500\\
\end{array} \right],
$$
i.e.,
$$
[R^{(1)}]^{\dag} =\left[ \begin{array}{c|c}
[R_1^{(1)}]^{\dag} & [R_2^{(1)}]^{\dag} \end{array} \right].
$$

Thus, the Moore\textendash Penrose inverse of $R^{(1)}$ can be calculated by finding the Moore\textendash  Penrose inverse of each block component.  $\hfill \Box$
\end{example}

\subsection{Application of EBOS to solve $YCC^* = AC^*$}
\label{ss:ebosycctact}

We will begin by describing the solution procedure for the system $YCC^* = AC^*$.  In other words, we wish to find $Y_+$  which satisfies $Y_+CC^* = AC^*$ and hence minimizes $\|A - YC\|_F$.  The general solution to this problem is given by
\begin{eqnarray}\label{eq:4.1}
 Y = Y_0 + K[I - (CC^*)^{\dag}]
\end{eqnarray}
where $Y_0 = AC^*[CC^*]^{\dag}$ is the minimum norm solution and $K \in {\mathbb C}^{m \times h}$ is arbitrary.  Our procedure for finding the proposed solution $Y_+$ is as follows.

Define $C^{(0)} = C$.  Suppose that $C^{(r-1)}$, for $r=1,\ldots, q-1$, has already been determined by
\begin{eqnarray}\label{eq:4.2}
C^{(r-1)} = \left[ \begin{array}{c}
C_1^{(0)} \\
\vdots \\
\rule{0cm}{0.5cm} C_{r}^{(r-1)} \\
\hline
\rule{0cm}{0.5cm} C_{r+1,c}^{(r-1)} \end{array} \right]\quad \text{where} \quad  C_{r+1, c}^{(r-1)}=\left[ \begin{array}{c}
C_{r+1}^{(r-1)} \\
\vdots \\
\rule{0cm}{0.5cm} C_q^{(r-1)} \end{array} \right].
\end{eqnarray}
To determine $C^{(r)}$ we define
\begin{eqnarray}\label{eq:4.3}
C_{r+1,c}^{(r)} = C_{r+1,c}^{(r-1)}\left( I -[ C_{r}^{(r-1)}]^{\dag} C_{r}^{(r-1)}\right)
\end{eqnarray}
and write
\begin{eqnarray}\label{eq:4.4}
C^{(r)} = \left[ \begin{array}{c}
C_1^{(0)} \\
\vdots \\
\rule{0cm}{0.5cm} C_{r}^{(r-1)} \\
\hline
\rule{0cm}{0.5cm} C_{r+1,c}^{(r)} \end{array} \right]\quad \text{where} \quad  C_{r+1,c}^{(r)} =\left[ \begin{array}{c}
C_{r+1}^{(r)} \\
\vdots \\
\rule{0cm}{0.5cm} C_q^{(r)} \end{array} \right].
\end{eqnarray}
To define $C^{(r+1)}$, for $r=0,\ldots, q-2$, the procedure in (\ref{eq:4.2})\textendash (\ref{eq:4.4}) is updated.
Further, let
\begin{eqnarray}\label{eq:4.5}
F_r = \left[ \begin{array}{cccc|c}
I & \bigzero & \cdots & \bigzero & \bigzero \\
\bigzero & I & \cdots & \bigzero & \bigzero \\
\vdots & \vdots & \ddots & \vdots & \vdots \\
\bigzero & \bigzero & \cdots & I & \bigzero \\ \hline
\rule{0cm}{0.5cm}\bigzero & \bigzero & \cdots & - C_{r+1,c}^{(r-1)} [ C_r^{(r-1)}]^{\dag} & I  \end{array} \right]
\end{eqnarray}
and $F = F_{q-1} \cdots F_1$.

\begin{theorem}\label{t:4.2}
For, $r=1,\ldots,q$, let $D_{C,rr}=C_r^{(r-1)} [C_r^{(r-1)}]^* $. Define
\begin{eqnarray}\label{eq:4.6}
D_C =\diag [ D_{C,11}, \dots, D_{C,qq}] \hspace{2mm} \mbox{and} \hspace{2mm} D^\dag_C =\diag [ D^\dag_{C,11}, \dots, D^\dag_{C,qq}].
\end{eqnarray}
 Then the matrix
\begin{eqnarray}\label{eq:4.7}
Y_+ = AC^* F^* D_C^\dag F
\end{eqnarray}
satisfies $Y_+CC^* = AC^*$ and hence also satisfies $\displaystyle\|A - Y_+C\|_F = \displaystyle\min _{Y} \|A - YC\|$.
\end{theorem}

\begin{remark}\label{r:1}
Matrix  $Y_+$ is not necessarily equal to $Y_0$ where
                             $\|Y_0\|_F = min \|Y\|_F$ subject to $YCC^* = AC^*$.
The detailed notes 1-5 at the end of the following proof explain this observation.
\end{remark}

\textbf{Proof.}  The first step is to define $C^{(0)} = C$ and explain the progressive block vector reduction
$$
C = C^{(0)} =\left[ \begin{array}{c}
C_1^{(0)} \\
\rule{0cm}{0.5cm} C_2^{(0)} \\
\vdots \\
\rule{0cm}{0.5cm} C_{q-1}^{(0)} \\
\rule{0cm}{0.5cm} C_q^{(0)} \end{array} \right] \quad \longrightarrow \quad
C^{(q-1)}= \left[ \begin{array}{c}
C_1^{(0)} \\
\rule{0cm}{0.5cm} C_2^{(1)} \\
\vdots \\
\rule{0cm}{0.5cm} C_{q-1}^{(q-2)} \\
\rule{0cm}{0.5cm} C_q^{(q-1)} \end{array} \right]
$$
to a block vector of mutually orthogonal blocks satisfying $C_r^{(r-1)} [C_s^{(s-1)}]^* = \bigzero$ for all $r \in \{1,\ldots,q\}$ and $s \in \{1,\ldots,q-1\}$ with $r > s$.  Thus the range $[C_s^{(s-1)}]^*({\mathbb C}^n)$ of $[C_s^{(s-1)}]^*$ is a subspace of the null space $[C_r^{(r-1)}]^{-1}( \{ \bfzero \})$ of $C_r^{(r-1)}$ for all $r > s$. That is $[C_s^{(s-1)}]^*({\mathbb C}^n) \subseteq [C_r^{(r-1)}]^{-1}( \{ \bfzero \})$ for all $r > s$.  At Stage $1$ we write
$$
C^{(0)} = \left[ \begin{array}{c}
C_1^{(0)} \\ \hline
\rule{0cm}{0.5cm} C_2^{(0)} \\
\vdots \\
\rule{0cm}{0.5cm} C_q^{(0)} \end{array} \right] = \left[ \begin{array}{c}
C_1^{(0)} \\ \hline
\rule{0cm}{0.5cm} C_{2,c}^{(0)} \end{array} \right]
$$
and use the consolidated partition to define
$$
F_1 = \left[ \begin{array}{c|c}
I & \bigzero \\ \hline
\rule{0cm}{0.5cm} - C_{2,c}^{(0)} [C_1^{(0)}]^{\dag} & I \end{array} \right].
$$
It follows that
$$
C^{(1)} = F_1C^{(0)} = \left[ \begin{array}{c}
C_1^{(0)} \\ \hline
\rule{0cm}{0.5cm} - C_{2,c}^{(0)} [C_1^{(0)}]^{\dag}C_1^{(0)} + C_{2,c}^{(0)} \end{array} \right] =  \left[ \begin{array}{c}
C_1^{(0)} \\ \hline
\rule{0cm}{0.5cm} C_{2,c}^{(1)} \end{array} \right] =  \left[ \begin{array}{c}
C_1^{(0)} \\ \hline
\rule{0cm}{0.5cm} C_2^{(1)} \\
\rule{0cm}{0.5cm} C_3^{(1)} \\
\vdots \\
\rule{0cm}{0.5cm} C_{q}^{(1)} \end{array} \right]
$$
where we have defined $C_{2,c}^{(1)} = C_{2,c}^{(0)} (I - [C_1^{(0)}]^{\dag}C_1^{(0)})$.  At Stage $2$ we write
$$
C^{(1)} = \left[ \begin{array}{c}
C_1^{(0)} \\
\rule{0cm}{0.5cm} C_2^{(1)} \\ \hline
\rule{0cm}{0.5cm} C_3^{(1)} \\
\vdots \\
\rule{0cm}{0.5cm} C_{q}^{(1)} \end{array} \right] = \left[ \begin{array}{c}
C_1^{(0)} \\
\rule{0cm}{0.5cm} C_2^{(1)} \\ \hline
\rule{0cm}{0.5cm} C_{3,c}^{(1)} \end{array} \right]
$$
and use the consolidated partition to define
$$
F_2 = \left[ \begin{array}{cc|c}
I & \bigzero & \bigzero \\
\bigzero & I & \bigzero \\ \hline
\rule{0cm}{0.5cm} \bigzero & - C_{3,c}^{(1)} [C_2^{(1)}]^{\dag} & I \end{array} \right].
$$
It follows that
$$
C^{(2)} = F_2C^{(1)} = \left[ \begin{array}{c}
C_1^{(0)} \\
\rule{0cm}{0.5cm} C_2^{(1)} \\ \hline
\rule{0cm}{0.5cm} - C_{3,c}^{(1)} [C_2^{(1)}]^{\dag}C_2^{(1)} + C_{3,c}^{(1)} \end{array} \right] = \left[ \begin{array}{c}
C_1^{(0)} \\
\rule{0cm}{0.5cm} C_2^{(1)} \\ \hline
\rule{0cm}{0.5cm} C_{3,c}^{(2)} \end{array} \right] = \left[ \begin{array}{c}
C_1^{(0)} \\
\rule{0cm}{0.5cm} C_2^{(1)} \\ \hline
\rule{0cm}{0.5cm} C_3^{(2)} \\
\rule{0cm}{0.5cm} C_4^{(2)} \\
\vdots \\
\rule{0cm}{0.5cm} C_q^{(2)} \end{array} \right]
$$
where we have defined $C_{3,c}^{(2)} = C_{3,c}^{(1)} (I - [C_2^{(1)}]^{\dag}C_2^{(1)})$.  In general, at Stage $r$, { the matrices  $C^{(r-1)}$  and $F_r$ are defined respectively by the expressions in (\ref{eq:4.2}) and (\ref{eq:4.5}).}   It follows that
$$
C^{(r)} = F_rC^{(r-1)} = \left[ \begin{array}{c}
C_1^{(0)} \\
\vdots \\
\rule{0cm}{0.5cm} C_r^{(r-1)} \\ \hline
\rule{0cm}{0.5cm}  C_{r+1,c}^{(r-1)}( I -[ C_r^{(r-1)}]^{\dag} C_r^{(r-1)}) \end{array} \right] = \left[ \begin{array}{c}
C_1^{(0)} \\
\vdots \\
\rule{0cm}{0.5cm} C_r^{(r-1)} \\ \hline
\rule{0cm}{0.5cm}  C_{r+1,c}^{(r)} \end{array} \right] = \left[ \begin{array}{c}
C_1^{(0)} \\
\vdots \\
\rule{0cm}{0.5cm} C_r^{(r-1)} \\ \hline
\rule{0cm}{0.5cm}  C_{r+1}^{(r)} \\
\vdots \\
\rule{0cm}{0.5cm} C_q^{(r+1)} \end{array} \right]
$$
where we have defined $C_{r+1,c}^{(r)} = C_{r+1,c}^{(r-1)}( I -[ C_r^{(r-1)}]^{\dag} C_r^{(r-1)})$.  The procedure continues until we obtain
$$
C^{(q-1)} =  \left[ \begin{array}{c}
C_1^{(0)} \\
C_2^{(1)} \\
\vdots \\
\rule{0cm}{0.5cm} C_q^{(q-1)} \end{array} \right].
$$
We can describe the final transformed matrix $C^{(q-1)}$ by the definitive relationships
$$
C_{s+1}^{(s)} = C_{s+1}^{(s-1)} [ I - \{ C_s^{(s-1)} \}^{\dag} C_s^{(s-1)} ]
$$
for each $s=1,\ldots,q-1$.  An important consequence of the reduction implemented by EBOS is that
\begin{eqnarray*}
C_{s+1,c}^{(s)} [C_s^{(s-1)}]^* & = & C_{s+1,c}^{(s-1)} (I - [C_s^{(s-1)}]^{\dag} C_s^{(s-1)}) [C_s^{(s-1)}]^* \\
& = & C_{s+1,c}^{(s-1)} (I -  [C_s^{(s-1)}]^* \{ C_s^{(s-1)} [C_s^{(s-1)}]^* \}^{\dag} C_s^{(s-1)}) [C_s^{(s-1)}]^* \\
& = & C_{s+1,c}^{(s-1)} (I -  [C_s^{(s-1)}]^* \{ [C_s^{(s-1)}]^* \}^{\dag} ) [C_s^{(s-1)}]^* \\
& = & \bigzero
\end{eqnarray*}
for each $s = 1,\ldots,q-1$.  Therefore $C_r^{(s-1)} [C_s^{(s-1)}]^*= \bigzero$ for each $s = 1,\ldots,q-1$ and all $r = s+1,\ldots,q$.  It follows that $C_{s+2,c}^{(s-1)} [C_{s+1}^{(s)}]^*= \bigzero$ for each $s = 1,\ldots, q-2$ and hence we deduce that
$$
C_{s+2,c}^{(s+1)} [C_s^{(s-1)}]^* = C_{s+2,c}^{(s)} (I - [C_{s+1}^{(s)}]^{\dag} C_{s+1}^{(s)}) [C_s^{(s-1)}]^* = \bigzero
$$
from which it follows that $C_r^{(s+1)}[C_s^{(s-1)}]^* = \bigzero$ for $r > s+1$ and in particular that $C_{s+2}^{(s+1)} [C_s^{(s-1)}]^* = \bigzero$.  By continuing this argument we can show that $C_r^{(r-1)} [C_s^{(s-1)}]^* = \bigzero$ for all $r > s$.  It follows that $C_s^{(s-1)} [C_r^{(r-1)}]^* = \bigzero$ for all $r > s$.  Therefore
$$
C^{(q-1)} [C^{(q-1)}]^* = \left[ \begin{array}{cccc}
C_1^{(0)} [C_1^{(0)}]^* & \bigzero &  \cdots & \bigzero \\
\bigzero & C_2^{(1)} [C_2^{(1)}]^* &  \cdots & \bigzero \\
\vdots & \vdots & \ddots & \vdots \\
\bigzero & \bigzero &  \cdots & C_q^{(q-1)} [C_q^{(q-1)}]^* \end{array} \right] = D_C.
$$
If we define $F = F_{q-1} \cdots F_1$ then $FC = C^{(q-1)}$ and
\begin{eqnarray}\label{eq:4.8}
FCC^*F^* = D_C.
\end{eqnarray}
If we also define $Z = YF^{-1}$ then
$$
YCC^* = AC^* \iff YF^{-1} \cdot FCC^*F^* = AC^*F^* \iff Z \cdot D_C = A[C^{(q-1)}]^*.
$$
Let ${\mathcal Y} = \{ Y \mid YCC^* = AC^* \}$ and ${\mathcal Z} = \{ Z \mid ZD_C = A[C^{(q-1)}]^*\}$ denote the respective solution sets.  We make the following observations.
\begin{enumerate}
\item $Y \in {\mathcal Y} \iff Z \in {\mathcal Z}$.
\item $Z_0 = A[C^{(q-1)}]^*D_C^{\dag} \Longrightarrow Z_0D_C = A[C^{(q-1)}]^* \iff Z_0F \cdot CC^* = AC^* \Longrightarrow Z_0F \in {\mathcal Y}$.
\item $Y_0 = AC^*(CC^*)^{\dag} \Longrightarrow Y_0CC^* = AC^* \iff Y_0F^{-1} D_C = A[C^{(q-1)}]^* \Longrightarrow Y_0F^{-1} \in {\mathcal Z}$.
\item $\|Z_0\| = \min_{Z \in {\mathcal Z}} \|Z\| \not \implies \|Z_0F\| = \min_{Y \in {\mathbb Y}} \|Y\| = \|Y_0\|$.
\item $\|Y_0\| = \min_{Y \in {\mathbb Y}} \|Y\| \not \implies \|Y_0F^{-1}\| = \min_{Z \in {\mathbb Z}} \|Z\| = \|Z_0\|$.
\end{enumerate}
 { In practice we will calculate $Y_+ = Z_0F = A [C^{q-1)}]^* [D_C]^{\dag}F = A [C^{(q-1)}]^{\dag}F$.  Note that in general $C^{\dag} \neq [C^{(q-1)}]^{\dag}F$.} It follows from the observations above that $Y_+CC^* = AC^*$ and hence that
$$
\| Z_0FC - A \| = \| Y_+C - A \| = \min_{Y \in {\mathcal S}} \| YC - A \| = \|Y_0C - A\|.
$$
Therefore $Z_0FC = Y_+C$ is an optimal reconstruction of $A$.  Note that in general $Y_+ \neq Y_0$. $\hfill \Box$

\subsection{Application of EBOS to solve $B^*BX = B^*Y$}\label{ss:bxy}

We wish to find $X_+$ which satisfies  $B^*BX_+ = B^*Y$ and hence minimizes $\|B^*Y - B^*BX\|_F$.  The general solution to the equation $B^*BX = B^*Y$ is given by
\begin{eqnarray}\label{eq:4.9}
X = X_0 + [I - (B^*B)^{\dag}]H
\end{eqnarray}
where $X_0 = (B^*B)^{\dag}B^*Y$ is the minimum norm solution and $H \in {\mathbb C}^{g \times h}$ is arbitrary.  Our procedure for finding the proposed solution {$X_+$} is as follows. 

Define $B^{(0)} = B$.  The transpose of the system $B^*BX = B^*Y \iff X^*B^*B = Y^*B$ has the same form as the system $YCC^* = AC^*$.  Thus we can implement an analogous reduction scheme.   Consequently we will restrict our remarks to a few key points.  We simplify the block row vector $B$ using right multiplication by a non-singular elementary block matrix to implement the desired column operations.  At Stage $r$ we obtain
\begin{eqnarray*}
&&\hspace*{-10mm}B^{(r-1)} = \left[ \begin{array}{ccc|c}
B_1^{(0)} & \cdots & B_r^{(r-1)} & B_{r+1,c}^{(r-1)} \end{array} \right]
\end{eqnarray*}
where $B_{r+1,c}^{(r-1)} = [ B_{r+1}^{(r-1)} \cdots B_p^{(r-1)}]$ and we use the consolidated partitions to define
$$
E_r = \left[ \begin{array}{ccccc}
I & \bigzero & \cdots & \bigzero & \bigzero \\
\bigzero & I & \cdots & \bigzero & \bigzero \\
\vdots & \vdots & \ddots & \vdots & \vdots \\
\bigzero & \bigzero & \cdots & I & -\{B_r^{(r-1)}\}^{\dag} B_{r+1,c}^{(r-1)} \\
\bigzero & \bigzero & \cdots & \bigzero & I \end{array} \right].
$$
Now we define $B^{(r)} = B^{(r-1)} E_r$.  The procedure terminates with
$$
B^{(p-1)} = [B_1^{(0)} B_2^{(1)} \cdots B_p^{(p-1)}].
$$
The final transformed matrix satisfies the definitive relationshipseq:4.10
$$
B_{s+1}^{(r)} = [ I - B_s^{(s-1)} \{ B_s^{(s-1)} \}^{\dag}] B_s^{(r)}
$$
for each $s=1,\ldots,p-1$ and each $r = s+1,\ldots,p-1$.  Similar arguments to those used previously for the block row reduction of the matrix $C$ can now be used to show that $[B_s^{(s-1)}]^*B_r^{(r-1)} = \bigzero$ and $[B_r^{(r-1)}]^*B_s^{(s-1)} = \bigzero$ for all $r > s$.  Therefore
$$
[B^{(p-1)}]^*B^{(p-1)} = \mbox{diag}([B_1^{(0)}]^*B_1^{(0)}, [B_2^{(1)}]^*B_2^{(1)}, \cdots, [B_p^{(p-1)}]^*B_p^{(p-1)})  = D_B
$$
is a block diagonal matrix.  If we define $E = E_1 \cdots E_{p-1}$ then $BE = B^{(p-1)}$ and $E^*B^*BE = D_B$.  If we also define $W = E^{-1}X$ then
$$
B^*BX = B^*Y \iff  E^*B^*BE \cdot E^{-1}X = E^*B^*Y \iff D_B \cdot W = [B^{(p-1)}]^*Y.
$$
We have the following theorem.

\begin{theorem}\label{t:4.3}
 {Let $D_{B,rr} = [B_r^{(r-1)}]^*B_r^{(r-1)}$ for each $r = 1, \ldots,p$ and define
 \begin{eqnarray}\label{eq:4.10}
 D_B = \mbox{diag}[ D_{B,11}, \ldots, D_{B,pp}] \quad \mbox{and}\quad {D_B}^{\dag} = \mbox{diag}[ {D_{B,11}}^{\dag}, \ldots, {D_{B,pp}}^{\dag}]
\end{eqnarray}
then the matrix
\begin{eqnarray}\label{eq:4.11}
X_+ = E{D_B}^{\dag}E^*B^*Y
\end{eqnarray}
satisfies $B^*BX_+ = B^*Y$ and hence minimizes $\| B^*Y - B^*BX \|_F$.  If $Y$ satisfies (\ref{eq:4.7}) then (\ref{eq:4.11}) solves the equation in (\ref{eq:1.3}).}  $\hfill \Box$
\end{theorem}

\textbf{Proof.} The transposed system $X^*B^*B = Y^*B \iff B^*BX = B^*Y$ has the same general form as the system $YCC^* = AC^*$.  Therefore the result follows by applying Theorem~\ref{t:4.2} to the transformed system.  A direct proof along similar lines to the proof of Theorem~\ref{t:4.2} is left to the reader.  $\hfill \Box$

Let ${\mathcal X} = \{ X \mid B^*BX = B^*Y\}$ and ${\mathcal W} = \{ W \mid D_BW = [B^{(p-1)}]^*Y \}$ denote the respective solution sets.  We make the following observations.
\begin{enumerate}
\item $X \in {\mathcal X} \iff W \in {\mathcal W}$.
\item $W_0 = D_B^{\dag}[B^{(p-1)}]^*Y \Rightarrow D_BW_0 = [B^{(p-1)}]^*Y \Leftrightarrow B^*B  \cdot EW_0 = BY^* \Rightarrow EW_0 \in {\mathcal X}$.
\item $X_0 = (B^*B)^{\dag}B^*Y \Rightarrow B^*BX_0 = B^*Y \Leftrightarrow D_B E^{-1}X_0 = [B^{(p-1)}]^*Y \Rightarrow E^{-1}X_0 \in {\mathcal W}$.
\item $\|W_0\| = \min_{W \in {\mathcal W}} \|W\| \Rightarrow \|EW_0\| = \min_{X \in {\mathbb X}} \|X\| = \|X_0\|$.
\item $\|X_0\| = \min_{X \in {\mathbb X}} \|X\| \Rightarrow \|E^{-1}X_0\| = \min_{W \in {\mathbb W}} \|W\| = \|W_0\|$.
\end{enumerate}
In practice we will define $W_0 = [D_B]^{\dag} [B^{(p-1)}]^*Y =  [B^{(p-1)}]^{\dag}Y$ and define $X_+ = EW_0$.  Note that in general $B^{\dag} \neq E[B^{(p-1)}]^{\dag}$.  It follows from the observations above that $B^*BX_+ = B^*Y$ and hence that
$$
\| BEW_0 - Y \| = \| BX_+ - Y \| = \min_{X \in {\mathcal X}} \| BX - Y \| = \| BX_0 - Y \|.
$$
Therefore $BEW_0 = BX_+$ is an optimal reconstruction of $Y$.

\subsection{EBOS algorithm}\label{ss:ebosa}
Here, we represent an algorithm for calculating $X_+$ as follows.
    
\begin{algorithm}[H]
Calculate $X_+$
\begin{algorithmic}
\REQUIRE $A, B, C, p, q, m, n, g_1,\ldots,g_p, h_1,\ldots, h_q$
\ENSURE $X_+$
\FOR{$r=1,\ldots,q-1$}
\STATE Calculate $C^{(r-1)}, C^{(r)}_{r+1,c}, C^{(r)}, F_r, F$
\ENDFOR
\FOR{$r=1,\ldots, q$}
\STATE Calculate $D_{C,rr}$
\ENDFOR
\STATE Calculate $D_C, D^{\dag}_C, Y_+$
\FOR{$r=1,\ldots,p-1$}
\STATE Calculate $B^{(r-1)}, B^{(r)}_{r+1,c}, B^{(r)}, E_r, E$
\ENDFOR
\FOR{$r=1,\ldots, p$}
\STATE Calculate $D_{B,rr}$
\ENDFOR
\STATE Calculate $D_B, D^{\dag}_B, X_+$

\end{algorithmic}
\end{algorithm}

\subsection{Analysis of numerical load associated with EBOS}
We calculate
$$
Y_+ = A [C^{(q-1)}]^* D_C F = A [ \{C_1^{(0)}\}^{\dag}, \cdots, \{C_q^{(q-1)}\}^{\dag} ] F.
$$
We assume $h_r = h/q$ for each $r=1,2,\ldots,q$.  

$\bullet$ \textbf{Number of flops to calculate \boldmath $C^{(q-1)}$ and $F$:}  Assume that $C^{(r-1)}$ is known.  To find $C^{(r)}$ we need to calculate
$$
[C_r^{(r-1)}]^{\dag} = [C_r^{(r-1)}]^* \{ C_r^{(r-1)} [C_r^{(r-1)}]^* \}^{\dag} \hspace*{2mm} \mbox{and} \hspace*{2mm} C_{r+1,c}^{(r)} = C_{r+1,c}^{(r-1)}( I -[ C_r^{(r-1)}]^{\dag} C_r^{(r-1)})
$$
for each $r = 1,\ldots, q-1$.  We have the following calculations.
\begin{enumerate}
\item $C_r^{(r-1)} [C_r^{(r-1)}]^*$ requires $2(h/q) \times n \times (h/q)$ flops;
\item $\{ [C_r^{(r-1)}]^*C_r^{(r-1)} \}^{\dag}$ requires $21(h/q)^3$ flops;
\item $[C_r^{(r-1)}]^{\dag} = [C_r^{(r-1)}]^*\{C_r^{(r-1)} [C_r^{(r-1)}]^* \}^{\dag}$ takes $2 \times n \times (h/q) \times (h/q)$ flops.
\item $[C_r^{(r-1)}]^{\dag} C_r^{(r-1)}$ requires $2 \times n \times (h/q) \times n$ flops;
\item $I -[ C_r^{(r-1)}]^{\dag} C_r^{(r-1)}$ involves only diagonal elements and requires $n$ flops; and
\item $C_{r+1,c}^{(r)} = C_{r+1,c}^{(r-1)}( I -[ C_r^{(r-1)}]^{\dag} C_r^{(r-1)})$ requires $2(q-r)h/q \times n \times n$ flops.
\end{enumerate}
Therefore stage $r$ requires $21h^3/q^3 + 4nh^2/q^2 + 2(q-r)n^2h/q + n + 2n^2h/q$ flops.  The total number of flops for all stages $r=1,\ldots,q-1$ to calculate $C^{(q-1)}$ is
$$
{\mathcal N}_1 = [21h^3/q^3 + 4nh^2/q^2 + 2n^2h/q + n^2h + n ] (q-1)
$$
We must also calculate $F = F_{q-1} \cdots F_1$.  We assume $Q_{r-1} = F_{r-1} \cdots F_1$ is known.  We can use the natural block structure to write $F_r = [F_{r,(i,j)}]_{i,j = 1}^q, Q_r = [Q_{r,(i,j)}]_{i,j = 1}^q$ where $F_{r,(i,j)}, Q_{r,(i,j)} \in {\mathbb C}^{(h/q) \times (h/q)}$ for each $r =1,\ldots,q-1$ and each $i,j = 1,\ldots,q$.  If we define
$$
F_r(r+1:q,r+1:q) = \left[ \begin{array}{cc}
I_{h/q} & \bigzero \\
-C_{r+1,c}^{(r-1)} [C_r^{(r-1)}]^{\dag} & I_{(q-r)h/q} \end{array} \right] \in {\mathbb C}^{(q-r)h/q \times (q-r+1)h/q}
$$
and
$$
Q_{r-1}(r:q,1:r) = \left[ \begin{array}{ccc}
Q_{r-1,(r,1)} & \cdots & Q_{r-1,(r,r)} \\
\vdots & \vdots & \vdots \\
Q_{r-1,(q,1)} & \cdots & Q_{r-1,(q,r)} \end{array} \right] \in {\mathbb C}^{(q-r+1)h/q \times rh/q}
$$
then to calculate $Q_r = F_rQ_{r-1}$ we need only calculate
$$
Q_r(r+1:q,1:r+1) = F_r(r+1:q,r+1:q)\,Q_{r-1}(r:q,1:r)
$$
for each $r=2,\ldots,q-1$.  This calculation requires $2(h/q)^3r(q-r)(q-r+1)$ flops.  The total number of flops for all stages is
\begin{eqnarray*}
{\mathcal N}_2 & = & 2(h/q)^3 \sum_{r=2}^{q-1} r(q-r)(q-r+1) \\
& = & (h/q)^3\left[ q^4/6 +q^3/3 -13q^2/6 + 5q/3 \right].
\end{eqnarray*}

$\bullet$ \textbf{Number of flops to calculate \boldmath $Y_+ = A [C^{(q-1)}]^*{D_C}^{\dag}F$:}  We have already calculated $\{C_r^{(r-1)} \}^{\dag}$ and so all matrices in the product are known.  The respective dimensions are $m \times n$, $n \times h$ and $h \times h$.  Therefore the number of flops required by the product is
$$
{\mathcal N}_3 = 2mnh + \min \{ 2mh^2, 2nh^2 \}.
$$

$\bullet$ \textbf{Number of flops used by EBOS to calculate \boldmath $Y_+$:} The total number of flops is given by ${\mathcal N} = {\mathcal N}_1 + {\mathcal N}_2 + {\mathcal N}_3$.

$\bullet$ \textbf{Number of flops used by the direct method to calculate \boldmath $Y_0 = A C^{\dag}$:}  We have $C \in {\mathbb C}^{h \times n}$ and so we calculate $C^{\dag} = C^*[CC^*]^{\dag}$.  Calculating $CC^*$ requires $2h^2n$ flops and  calculating $[CC^*]^{\dag}$ takes $21h^3$ flops.  Calculating $C^*[CC^*]^{\dag}$ takes another $2nh^2$ flops and finally $AC^{\dag}$ takes a further $2mnh$ flops.  Therefore the total number of flops for the direct method\footnote{This is not necessarily a true indication of the time taken by M{\sc atlab} to compute $Y_0 = AC^{\dag}$.} is ${\mathcal F} = 21h^3 + 4nh^2 + 2mnh$.

$\bullet$ \textbf{A typical comparison:}  We wish to compare typical values for ${\mathcal F}$ and ${\mathcal N}$.  Suppose $m = n$ and $h = \epsilon n$ where $\epsilon \in (0,1]$ and $q=5$.  We assume $n$ is large and so we only include the highest order terms in $n^3$.  For EBOS ${\mathcal N}_1 \approx  (21 \epsilon^3/q^3 + 4\epsilon^2/q^2 + 2\epsilon/q + \epsilon)n^3$, ${\mathcal N}_2 \approx  \epsilon^3[q/6 +1/3 - 13/(6q) + 5/(3q^2)]n^3$ and ${\mathcal N}_3 = (2 \epsilon + 2 \epsilon^2)n^3$.  For the direct method ${\mathcal F} = (21 \epsilon^3 + 4 \epsilon^2 + 2 \epsilon)n^3$.

\begin{table}[htb]
\begin{center}
\vspace{0.25cm}
\begin{tabular}{|c|c|c|c|c|c|c|c|} \hline\hline
$\epsilon$ & $q$ & ${\mathcal N}_1/n^3$ & ${\mathcal N}_2 /n^3$ & ${\mathcal N}_3 /n^3$ & ${\mathcal N}/n^3$ & ${\mathcal F}/n^3$ & ${\mathcal N}/{\mathcal F}$ \\ \hline \hline
1 & 2 & 5.62 & 1.67 & 4 & 11.29 & 27 & 0.4182 \\ \hline
1 & 3 & 2.89 & 2.80 & 4 & 9.69 & 27 & 0.3587 \\ \hline
1 & 5 & 1.73 & 4.97 & 4 & 10.69 & 27 & 0.3961 \\ \hline
1 & 10 & 1.26 & 10.13 & 4 & 15.39 & 27 & 0.5702 \\ \hline \hline
0.75 & 2 & 3.17 & 0.70 & 2.63 & 6.50 & 12.61 & 0.5153 \\ \hline
0.75 & 3 & 1.83 & 1.18 & 2.63 & 5.63 & 12.61 & 0.4467 \\ \hline
0.75 & 5 & 1.21 & 2.10 & 2.63 & 5.93 & 12.61 & 0.4704 \\ \hline
0.75 & 10 & 0.93 & 4.28 & 2.63 & 7.83 & 12.61 & 0.6211 \\ \hline \hline
0.5 & 2 & 1.58 & 0.21 & 1.5 & 3.29 & 4.63 & 0.7106 \\ \hline
0.5 & 3 & 1.04 & 0.35 & 1.5 & 2.89 & 4.63 & 0.6251 \\ \hline
0.5 & 5 & 0.76 & 0.62 & 1.5 & 2.88 & 4.63 & 0.6231 \\ \hline
0.5 & 10 & 0.61 & 1.27 & 1.5 & 3.38 & 4.63 & 0.7307 \\ \hline \hline
0.25 & 2 & 0.60 & 0.03 & 0.63 & 1.25 & 1.08 & {\color{red}1.1636} \\ \hline
0.25 & 3 & 0.46 & 0.04 & 0.63 & 1.13 & 1.08 & {\color{red}1.0437} \\ \hline
0.25 & 5 & 0.36 & 0.08 & 0.63 & 1.07 & 1.08 & 0.9880 \\ \hline
0.25 & 10 & 0.30 & 0.16 & 0.63 & 1.09 & 1.08 & {\color{red}1.0075} \\ \hline \hline
\end{tabular}
\end{center}
\vspace{0.25cm}
\caption{Number of flops for EBOS (${\mathcal N}$) and the direct method (${\mathcal F}$) to solve $YCC^* = AC^*$ when $m = n$ and $h = \epsilon n$.  Note that EBOS is less effective when $\epsilon$ is small and $q$ is large.}
\label{tab1}
\end{table}

Table~\ref{tab1} shows that if $m=n$ and $h = \epsilon n$ the advantage for EBOS is greatest when $\epsilon \in [0.75, 1]$.   This is not surprising because the direct calculation becomes much easier when $h$ is relatively small. 

\begin{example}\label{ex:2}
In this example {\color{red} we} illustrate the advantage of EBOS over the direct method (\ref{eq:1.4}) for large matrices of particular sizes. We used the Matlab code in Section 7 to compare the computation times for solution of the equation $YCC^* = AC^*$ with different randomly chosen matrices $A$ and $C$.  In the code we have used the notation $h_j = dh = h/q$ for each $j = 1,\ldots,q$.  The random choice of matrices means that the computation times vary from one experiment to the next but the variation is relatively small compared to the execution time. In the tests we used $m = n$ with $h = n \iff \epsilon = 1$ and various $n \in [1000, 15000]$ and $q \in [2, 10]$.   Selected typical results are given in Tables~\ref{tab2}, \ref{tab3} and \ref{tab4}.  We have not bothered to illustrate solution of the second equation $B^*BX = B^*Y$ because it is essentially the same procedure used for the first equation.

\begin{table}[htb]
\begin{center}
\vspace{0.25cm}
\begin{tabular}{|c|c|c|c|c|c|} \hline\hline
$m=n=h$ & $dh$ & q & $t_{\mbox{\scriptsize d}}$ (secs)& $t_{\mbox{\scriptsize e}}$ (secs) & $t_{\mbox{\scriptsize e}}/t_{\mbox{\scriptsize d}}$ \\ \hline\hline
3000 & 1500 & 2 & 7.2597 & 6.9437 & 0.9565  \\ \hline
3000 & 1000 & 3 & 7.2267 & 6.2511 & 0.8650  \\  \hline
3000 & 750 & 4 & 7.5824 & 6.0482 & 0.7977  \\ \hline
3000 & 600 & 5 & 6.8432 & 5.8089 & 0.8489 \\ \hline
3000 & 500 & 6 & 7.3092 & 6.2032 & 0.8487 \\ \hline
3000 & 375 & 8  & 7.4798 & 6.8610 & 0.9173 \\ \hline
3000 & 300 & 10 & 7.4431 & 7.6427 & {\color{red}1.0268} \\  \hline
\end{tabular}
\end{center}
\vspace{0.25cm}
\caption{Computation time in seconds for the direct method and EBOS with $m = n = h = 3000$.}
\label{tab2}
\end{table}

\begin{table}[htb]
\begin{center}
\vspace{0.25cm}
\begin{tabular}{|c|c|c|c|c|c|} \hline\hline
$m=n=h$ & $dh$ & q & $t_{\mbox{\scriptsize d}}$ (secs)& $t_{\mbox{\scriptsize e}}$ (secs) & $t_{\mbox{\scriptsize e}}/t_{\mbox{\scriptsize d}}$ \\ \hline\hline
9000 & 4500 & 2 & 215.9239 & 180.0193 & 0.8377  \\ \hline
9000 & 3000 & 3 & 211.9872 & 155.5800 & 0.7339  \\  \hline
9000 & 2250 & 4 & 222.0178 & 144.5274 & 0.6510  \\ \hline
9000 & 1800 & 5 & 218.2690 & 138.2073 & 0.6332  \\ \hline
9000 & 1500 & 6 & 217.9000 & 144.1707 & 0.6616  \\ \hline
9000 & 1125 & 8  & 220.7669 & 155.2825 & 0.7034  \\ \hline
9000 & 900 & 10 & 211.7547 & 173.5781 & 0.8197 \\  \hline
\end{tabular}
\end{center}
\vspace{0.25cm}
\caption{Computation time in seconds for the direct method and EBOS with $m = n = h = 9000$.}
\label{tab3}
\end{table}

\begin{table}[htb]
\begin{center}
\vspace{0.25cm}
\begin{tabular}{|c|c|c|c|c|c|} \hline\hline
$m=n=h$ & $dh$ & q & $t_{\mbox{\scriptsize d}}$ (secs)& $t_{\mbox{\scriptsize e}}$ (secs) & $t_{\mbox{\scriptsize e}}/t_{\mbox{\scriptsize d}}$ \\ \hline\hline
15000 & 7500 & 2 & 1058.4895 & 930.0785 & 0.8787  \\ \hline
15000 & 5000 & 3 & 1128.0494 & 827.3659 & 0.7334  \\  \hline
15000 & 3750 & 4 & 1051.1613 & 706.4415 & 0.6721  \\ \hline
15000 & 3000 & 5 & 1031.2660 & 691.3852 & 0.6704  \\ \hline
15000 & 2500 & 6 & 1055.1534 & 706.5885 & 0.6697 \\ \hline
15000 & 1875 & 8  & 1178.7413 & 843.8233 & 0.7159 \\ \hline
15000 & 1500 & 10 & 1054.4458 & 903.6170 & 0.8570 \\  \hline
\end{tabular}
\end{center}
\vspace{0.25cm}
\caption{Computation time in seconds for the direct method and EBOS with $m = n = h = 15000$.}
\label{tab4}
\end{table}

For optimal choices of $q$ the results in Tables~\ref{tab2} \ref{tab3} and \ref{tab4} suggest that EBOS may reduce computation time by as much as $40\%$ compared to the direct method.
\end{example}

\begin{example}\label{ex2}
Details of the EBOS computation for a particular choice of small matrices.  Let
$$
A = \left[ \begin{array}{cccccccc}
1   &  1  &   1  &   1   &  1   &  1   &  1  &   1 \\
1   &  0   &  1  &  0   &  1  &   0  &   0   &  0 \\
0   &  1   &  1   &  1   &  0   &  1  &   1  &   1 \\
1  &   0   &  0   &  0   &  1   &  1  &   0  &   1 \\
1   &  0   &  1   &  0   &  1   &  1  &   0   &  0 \\
0   &  1  &   0   &  0   &  0   &  1   &  0   &  1 \\
0   &  1  &  1   &  1   &  0   &  1   &  1   &  0 \\
1   &  1  &  0   &  1  &   0  &   0   &  0  &  0 \\
1   &  1  &   0  &   0  &   1  &   0  &   0   &  1 \\
1   &  0  &  0   &  0  &   0   &  0   &  1   &  1 \\
0   &  1   &  0  &   0  &  1   &  1   &  1  &   1 \\
1   &  1   &  1   &  1   &  0  &   0   &  1  &   0 \end{array} \right], \ B = [B_1 \mid B_2] = \left[ \begin{array}{cc|cccc}
1   &  1   &  0   &  0   &  0   &  0 \\
0   &  1   &  1   &  0   &  0  &   1 \\
0   &  1   &  0  &   0   &  1   &  0 \\
0  &   0   &  1  &   0  &   0  &   0 \\
0   &  0   &  1   &  1   &  0   &  1 \\
1  &   1  &   1  &   1  &   0   &  0 \\
0  &   0  &   0  &   1  &   0   &  0 \\
1   &  1  &  0  &   0  &   0   &  0 \\
0   &  1   &  0   &  1   &  1  &   0 \\
1   &  1  &   1  &   1  &   1  &   0 \\
0  &   0   &  0  &   0  &   0   &  1 \\
1   &  0   &  1   &  1   &  0   &  1 \end{array} \right],
$$
$$
C = \left[ \begin{array}{c}
C_1 \\ \hline
C_2 \\ \hline
C_3 \end{array} \right] = \left[ \begin{array}{cccccccc}
1  &   0   &  1   &  0   &  1  &   1   &  0   &  0 \\
0  &  1  &   0   &  1  &   1  &   0  &   0  &   0 \\ \hline
0   &  1   &  1  &   1  &   1   &  1  &   0  &   0 \\
0   &  1   &  0  &   0   &  1  &   1  &   0  &   0 \\ \hline
1   &  0   &  0   &  0   &  1   &  1  &   1   &  1 \\
0   &  1  &   1   &  0  &   0   &  1  &   0   &  0 \\
0   &  0   &  1   &  0   &  1   &  1   &  0  &   0 \end{array} \right].
$$
\textbf{ \textit{Solution of \boldmath $YCC^* = AC^*$}}
We calculate
$$
F_1 =  \left[ \begin{array}{c|c}
I & \bigzero \\ \hline
\rule{0cm}{0.5cm} - C_{2c}^{(0)} [C_1^{(0)}]^{\dag} & I \end{array} \right] = \left[ \begin{array}{rr|rrrrr}
1    &     0    &     0    &     0    &     0     &    0     &    0 \\
0   & 1    &     0    &     0   &     0    &     0     &    0 \\ \hline
- 6/11 &  - 9/11  &  1   &      0    &     0     &    0     &    0 \\
- 4/11  & - 6/11   &      0  &  1   &      0   &     0    &     0 \\
- 8/11  & - 1/11    &     0    &     0  &  1     &    0    &     0 \\
- 5/11 &  - 2/11    &     0     &    0    &     0   & 1    &     0 \\
- 8/11  & - 1/11    &     0    &     0     &    0     &    0  &  1 \end{array} \right]
$$
and $C^{(1)} = F_1C$ giving
$$
C^{(1)} = \left[ \begin{array}{rr|rrrrrr}
1    &    0  &  1    &     0  &  1  &  1    &     0    &     0 \\
0  &  1    &     0  &  1  &  1   &     0     &    0    &     0 \\ \hline
- 5/11   &  2/11  &  5/11  & 2/11 &  - 4/11  &  5/11   &     0    &     0 \\
- 4/11  &  5/11 &  - 4/11  & - 6/11  &  1/11  &  7/11    &     0    &     0 \\
3/11 &  - 1/11 &  - 8/11 & - 1/11  &  2/11  &  3/11  &  1  &  1 \\
- 5/11 &  9/11  &  6/11 &  - 2/11  & - 7/11  &  6/11    &     0    &     0 \\
- 8/11  & - 1/11  &  3/11  & - 1/11  &  2/11  &  3/11   &     0    &     0 \end{array} \right].
$$
Next we calculate
$$
F_2 = \left[ \begin{array}{c|c|c}
I & \bigzero & \bigzero \\ \hline
\bigzero & I & \bigzero \\ \hline
\rule{0cm}{0.5cm} \bigzero & - C_{3c}^{(1)} [C_2^{(1)}]^{\dag} & I \end{array} \right] = \left[ \begin{array}{rr|rrr|rr}
1     &    0     &    0     &    0    &     0   &      0    &     0 \\
0   & 1     &    0    &     0    &     0    &     0    &     0 \\ \hline
0     &    0  &  1   &      0    &    0    &     0     &    0 \\
0     &    0     &    0  &  1    &     0    &     0    &     0 \\
0    &     0    &     0     &    0  &  1    &     0    &     0 \\ \hline
0    &     0  & -26/25  & - 2/5  &  2/25  &  1    &     0 \\
0     &    0 &  -13/25  & - 1/5  &  1/25    &     0  &  1 \end{array} \right]
$$
and $C^{(2)} = F_2C^{(1)}$ giving
$$
C^{(2)} = \left[ \begin{array}{rr|rrr|rrr}
1    &    0  & 1    &     0  &  1   & 1    &     0     &    0 \\
0  &  1     &    0  &  1  &  1    &     0    &     0     &    0 \\ \hline
- 6/11  &  2/11  &  5/11  &  2/11 &  -4/11  &  5/11   &      0   &      0 \\
- 4/11  &  5/11  & - 4/11  & - 6/11  &  1/11  &  7/11    &     0    &     0 \\
3/11  & - 1/11  & - 8/11  & - 1/11  &  2/11  &  3/11  &  1  &  1 \\ \hline
7/25  &  11/25  &  4/25  & - 4/25 &  - 7/25  & - 4/25  &  2/25  &  2/25 \\
- 9/25  & - 7/25  &  2/25  & - 2/25 &  9/25  &  - 2/25  &  1/25 &   1/25 \end{array} \right].
$$
Finally we calculate $Y_+ = A[C^{(2)}]^* [D_C]^{\dag}F_2F_1 = A[C^{(2)}]^{\dag} F$ giving
$$
Y_+ = \left[ \begin{array}{rr|rrr|rr}
0   & 1 &  0  &  -1  & 1  &  1 &  0 \\
1  &  2  & -2  & -1  &  0  &  1  &  1 \\
-1 &  0   & 1 &  -1 &   1 &  1 &  0 \\
0.5 &  -0.5 &   0.5 &   0.5 &   0.5 &  -0.5 &   -0.5 \\
1 &   0 &  0 &  0 &  0 &   0  &  0 \\
-0.5 &  -0.5 &   0.5 &   0.5 &   0.5 &   0.5 &  -0.5 \\
-0.5 &  -0.5 &   1.5 &  -0.5 &   0.5 &   0.5 &  -0.5 \\
1 &  0   & 1 &   0  & 0  & 0 &  -2 \\
0.5 &   2.5 &  -2.5 &  -0.5 &   0.5 &   1.5 &   0.5 \\
0  &  1 &  -1 &  -1 &   1 &   1 &  0 \\
-1 &   1 & -1 &  0  &   1 &   1 &   1 \\
0.5 &   1.5 &  -0.5 &  -1.5 &   0.5 &   1.5 &  -0.5 \end{array} \right].
$$

\textbf{\textit{Solution of \boldmath{$B^*BX = B^*Y_+$}}}.  Define
$$
E_1 = \left[ \begin{array}{rr|rrrr}
1 &  0  & -0.4737  & -0.4737  &  0.2632  & -0.1579 \\
0  &  1 &  -0.1579 &  -0.1579  & -0.5789  & -0.0526 \\ \hline
0     &    0  &  1   &      0    &     0    &     0 \\
0    &     0     &    0  &  1   &      0     &    0 \\
0    &     0     &    0    &     0  &  1   &      0 \\
0    &     0    &     0     &    0    &     0  &  1  \end{array} \right]
$$
and calculate $B^{(1)} = BE_1$ giving
$$
B^{(1)} = \left[ \begin{array}{rr|rrrr}
1 &  1 &   -0.6316  & -0.6316 & -0.3158  & -0.2105 \\
0  & 1 &    0.8421 &  -0.1579  & -0.5789  &  0.9474 \\
0   & 1 &  -0.1579  & -0.1579  &  0.4211 &  -0.0526 \\
0     &    0  &  1     &    0     &    0    &     0 \\
0     &    0  &  1 &   1  &       0   & 1 \\
1 &   1 &   0.3684  &  0.3684 &  -0.3158 &  -0.2105 \\
0    &     0     &    0   & 1    &     0   &      0 \\
1 &   1 &  -0.6316 &  -0.6316  & -0.3158 &  -0.2105 \\
0  &  1 &  -0.1579  &  0.8421  &  0.4211 &  -0.0526 \\
1 &   1 &   0.3684  &  0.3684  &  0.6842  & -0.2105 \\
0     &    0    &     0    &     0    &     0  &  1 \\
1   &      0  &  0.5263  &  0.5263  &  0.2632 &  0.8421 \end{array} \right].
$$
Define $E = E_1$.  Now $X_+ = E[D_B]^{\dag}[B^{(1)}]^*Y = E[B^{(1)}]^{\dag}Y$ which gives
$$
X_+ = \left[ \begin{array}{rr|rrr|rr}
 -0.1987  &  0.1009  &  0.4826  & -0.5599  &  0.4527  &  0.3691  & -0.7823 \\
 0.2934  &  0.4700  & -0.2839  &  0.0647  & -0.0016  &  0.3123  &  0.1073 \\ \hline
 0.4621 &  -0.5284  &  0.0205 &   0.2981  & -0.0804 & -0.5726  & -0.0300 \\
 -0.0126  & -0.0095  &  0.1735  &  0.0994 &  -0.0268  &  0.3091 &  -0.1767 \\
 -0.5394  &  0.8454  & -0.8328  & -0.8770  &  0.7287 &  0.7161 &   0.4479 \\
 0.0110  &  1.3833 &  -1.0268  & -0.7744  &  0.3360   & 0.9795  &  0.6546 \end{array} \right].
$$
The partitions indicate the individual $X_{+,ij}$ for $i=1,2$ and $j=1,2,3$. $\hfill \Box$
\end{example}

\section{Conclusion}

{ In this paper, we have developed a method called} the elementary block operations scheme (EBOS)  for calculating the optimal approximation of an observed matrix $A\in\ct^{m\times n}$. We targeted the case when matrix $A$ was large. The approximant is represented by $\sum_{i=1}^p \sum_{j=1}^q B_iX_{ij}C_j$  where each $B_i \in \ct^{m \times g_i}$ and $C_j \in \ct^{h_j \times n}$ is known and matrix $ X_{ij}$ is unknown. Each $ X_{ij}$ should be determined from the solution of the problem of minimizing the Frobenius norm of the error. The problem is motivated by the task in optimal data processing. Our solutions are motivated by the  observation that, for large $A$, the known approaches would involve considerable computational burden needed for calculation of large pseudo-inverse matrices. This would make computation unwieldy and inefficient.

The idea of our approach is to reduce the solution procedure to finding the pseudo-inverses of a collection of much smaller matrices.  In EBOS this idea is realized by using a sequence of elementary block row operations to reduce the associated systems to block diagonal form.  As a result, the computational time is reduced by as much as $40\%$.  The proposed scheme can be easily implemented using the standard general purpose M{\sc atlab} algorithms.

\begin{appendix}
\section{ M{\sc atlab} code for comparison of the direct method and EBOS} \label{s:mc}

\begin{lstlisting}[style=Matlab-editor]
clear all \\

%solve Y *C*C' = A*C', size(A) = (m,n), size(C) = (q*dh,n). \\

%define the parameters \\

m = 18; n = 18; dh = 3; q = 6; \\
$\%$m = 3000; n=3000; dh = 750; q = 4; \\

%define the random matrices \\

tic \\
A = rand(m, q*dh)*rand(q*dh,n); \\
C = rand(q*dh,n); \\
trand = toc \\ }

%set C0 = C for checking

C0 = C;

%calculate Y by the direct method

tic
Ydirect = A*pinv(C);
tdirect = toc
edirect = A - Ydirect*C;
ndirect = max(max(abs(edirect)))

%calculate Y by EBOS

C = C0;
tic
F = eye(q*dh);
F(dh+1:q*dh,1:dh) = - C(dh+1:q*dh,:)*pinv(C(1:dh,:));
C(dh+1:q*dh,:) = C(dh+1:q*dh,:)*(eye(n) - pinv(C(1:dh,:))*C(1:dh,:));
D(1:dh,1:dh) = C(1:dh,1:n)*C(1:dh,1:n)';
Ddag = pinv(D(1:dh,1:dh));
for r=2:q-1
F(r*dh+1:q*dh,1:r*dh) = [- C(r*dh+1:q*dh,:)*pinv(C((r-1)*dh+1:r*dh,:)),eye((q-r)*dh)]*...
    F((r-1)*dh+1:q*dh,1:r*dh);
C(r*dh+1:q*dh,:) = C(r*dh+1:q*dh,:)*(eye(n) - pinv(C((r-1)*dh+1:r*dh,:))*C((r-1)*dh+1:r*dh,:));
D((r-1)*dh+1:r*dh,(r-1)*dh+1:r*dh) = C((r-1)*dh+1:r*dh,:)*C((r-1)*dh+1:r*dh,:)';
Ddag = blkdiag(Ddag,pinv(D((r-1)*dh+1:r*dh,(r-1)*dh+1:r*dh)));
end
D((q-1)*dh+1:q*dh,(q-1)*dh+1:q*dh) = C((q-1)*dh+1:q*dh,:)*C((q-1)*dh+1:q*dh,:)';
Ddag = blkdiag(Ddag,pinv(D((q-1)*dh+1:q*dh,(q-1)*dh+1:q*dh)));
ACp = A*C';
ACpDd = [];
for r=1:q
    ACpDd = [ACpDd ACp(:,(r-1)*dh+1:r*dh)*Ddag((r-1)*dh+1:r*dh,(r-1)*dh+1:r*dh)];
end
Yebos = [];
for r=1:q
Yebos = [Yebos ACpDd(1:m,(r-1)*dh+1:q*dh)*F((r-1)*dh+1:q*dh,(r-1)*dh+1:r*dh)];
end
tebos = toc
eebos = Yebos*C0 - A;
nebos = max(max(abs(eebos)))
\end{lstlisting}

\end{appendix}

\end{document}